\documentclass{amsart}
\usepackage{amssymb,amsmath}
\DeclareMathOperator{\sech}{sech}
\DeclareMathOperator{\csch}{csch}
\newcommand{\beq}{\begin{equation}}
\newcommand{\eeq}{\end{equation}}
\newcommand{\beqstar}{\begin{equation*}}
\newcommand{\eeqstar}{\end{equation*}}
\newcommand{\del}{\partial}
\newcommand{\ket}{\rangle}
\newcommand{\bra}{\langle}

\begin{document}

\title[Classification of Einstein Metrics]{Classification of Einstein metrics on $I \times S^3$}
\date{12/6/2004}
\author[C.T. Asplund]{Curtis T. Asplund}
\curraddr[C.T. Asplund]{Department of Physics, University of California, Santa Barbara, CA, 93106}
\email[C.T. Asplund]{casplund@physics.ucsb.edu}
\author[B. Krummel]{Brian Krummel} 
\author[E.D. Merrell]{Evan D. Merrell} 
\author[R. Rachel]{Robert T. Rachal}
\author[D. Yang]{DaGang Yang}

\address{Department of Mathematics, Tulane University, 6823 St. Charles Ave., New Orleans, LA 70118}



\thanks{We would like to thank Tulane University for hosting us, and the National Science
	Foundation for providing the Research Experiences for Undergraduates program which supported this research.}

\begin{abstract}
We present a complete classification of Einstein metrics on the space 
	$M = I \times S^3$, where $I$ is the interval $(0,l)$ or $(0,\infty)$ or their closures, 
	and we consider separate metric functions $f$ and $h$
	of $t \in I$ for the base and fiber of the Hopf fibration $S^1 \rightarrow  S^3 \rightarrow S^2$.
	All such metrics yielding smooth and complete manifolds are included and discussed. 
	The results are surprisingly rich, including many well-known examples 
	and several one-parameter families of metrics with a variety of geometries.
\end{abstract}
	
\maketitle
	
\section{Introduction}

	Einstein manifolds in four dimensions are important in both geometry and in physics, where, as 
gravitational instantons, they have applications in quantum gravity (see \cite{Gibb} and \cite{Page}). A well-studied 4-manifold
		is the cylinder on the 3-sphere $M = I \times S^3$, $I$ an interval of the real line.
		The Hopf fibration provides a natural way of decomposing the geometry of $S^3$ into two components, the base space
		$S^2$ and fibers $S^1$. We define real-valued functions $h(t)$ and $f(t)$ such that
		for each $t \in I$, these functions determine the metric on the submerged components in $S^3$. 
		We then pose and answer the question: what functions $f$ and $h$ yield Einstein manifolds $(M,g)$? 
 
	The first example of a compact, \emph{inhomogenous} Einstein 4-manifold (with positive scalar curvature) 
		was given in 1978 by D. Page \cite{Page}, a non-trivial fibration of $S^2$ by $S^2$. This 
		example appears in our classification in section \ref{Page's example}.
		Soon afterwards, B$\acute{\text{e}}$rard Bergery (see \cite{Berg} and \cite{Besse}) gave a generalization of
		Page's result to arbitrary even dimension.
		We concentrate on dimension four. The key to our approach is a change of independent variable
		that allows us to isolate $h$ and integrate it explicitly in terms of this new variable and a small 
		number of parameters. The bulk of the paper is an exploration of parameter space in search of solutions 
		that yield smooth and geodesically complete Einstein manifolds. In a number of cases, we are able to integrate 
		and obtain explicit formulas. In the non-integrable cases, we may still determine if complete non-singular
		manifolds are possible, allowing us to classify all possible solutions of interest.
	
	We begin in section 2 with a discussion of the global frame and metric we use for $M$ and 
		a calculation of its Ricci curvature
		tensor in terms of the functions $f$ and $h$. 
		Restricting to Einstein metrics, we obtain ODE's that admit explicit general solutions via the change of
		variable mentioned above.
		The classification and discussion of solutions occurs in sections 3 through 6, and we conclude in section 7.
		Along the way we encounter a rich collection of spaces,
		including $\mathbb{R}^4$, $S^4$, $\mathbb{C}P^2$, $TS^2$,
		$\mathbb{C}P^2 \# \overline{\mathbb{C}P^2}$ (Page's example), 
		several one-parameter families of Einstein metrics, and Einstein orbifolds.
	
\section{Preliminaries}

	We begin by constructing a global frame for $M = I \times S^3$. 
		Considering the quaternions $\mathbb{H} = \mathbb{C} \times \mathbb{C}$, a real
		vector space with basis $\mathbf{1} = (1,0)$, $\mathbf{i} = (i,0)$, $\mathbf{j} = (0,1)$, $\mathbf{k} = (0,i)$,
		and inner product given by $\bra p,q \ket = \frac{1}{2} (p^*q + q^*p)$, where $p^* = (\bar{a},-b)$ for $p = (a,b)$.
		The function
		$\phi :\mathbb{H} \rightarrow \mathbb{R}^4$ given by $\phi (x^1 + x^2 i, x^3 + x^4 i) = (x^1,x^2,x^3,x^4)$ 
		provides a natural smooth manifold structure on $\mathbb{H}$; it is just the structure of $\mathbb{R}^4$ obtained
		from $\phi^{-1}$.
		Define vector fields given at each point $q \in \mathbb{H}$ by  
		$\mathcal{X}_1|_q = \mathbf{i}q, \mathcal{X}_2|_q = \mathbf{j}q, \mathcal{X}_3|_q = \mathbf{k}q$, 
		where we multiply quaternions according to $(a,b)(c,d) = (ac - \bar{d}b, da + b\bar{c})$, 
		and make the canonical identification $T_q\mathbb{H} = \mathbb{H}$.
		Let $\mathcal{S} = \{q \in \mathbb{H} \mid \bra q,q \ket = 1 \}$, the three dimensional submanifold of unit 
		quaternions.
		Since for any ``imaginary" $p \in \mathbb{H}$ (i.e. $p^* = -p$), $\bra q,pq \ket = 0$, the set 
		$\{ \mathcal{X}_i|_p \}$ spans each tangent space $T_pM$ and $\{ \mathcal{X}_i \}$ is a global frame for 
		$\mathcal{S}$
		(see \cite{Lee} p.203).
		
	The Hopf fibration is associated with the partition of $\mathcal{S}$ by unit-speed curves (the Hopf fibers)
		defined
		for each $q = (a,b) \in \mathcal{S}$ by $c_q(\tau) = (ae^{i\tau}, be^{i\tau})$, $\tau \in \mathbb{R}$ 
		(see \cite{Walschap} p.29).
		Note that $c' = \mathbf{i}q$, thus $\mathcal{X}_1$ is tangent to the Hopf fiber at each point.
		
		Now $\phi$ restricts to a diffeomorphism $\phi : \mathcal{S} \rightarrow S^3$ that pushes forward to a map 
		$\phi_*: T\mathcal{S} \rightarrow TS^3$ between tangent bundles that yields a global frame for $S^3$:
		\begin{align*}
			X_1 &= -x^2\frac{\del}{\del x^1} + x^1\frac{\del}{\del x^2} - x^4\frac{\del}{\del x^3} + x^3\frac{\del}{\del 
				x^4}, \\
			X_2 &= -x^3\frac{\del}{\del x^1} + x^4\frac{\del}{\del x^2} + x^1\frac{\del}{\del x^3} - x^2\frac{\del}{\del
				x^4}, \\
			X_3 &= -x^4\frac{\del}{\del x^1} - x^3\frac{\del}{\del x^2} + x^2\frac{\del}{\del x^3} + x^1\frac{\del}{\del
				x^4},
		\end{align*}
		where $X_i = \phi_*(\mathcal{X}_i)$.
		For arbitrary smooth vector fields $V = V^i \del /\del x^i$, $W = W^j \del /\del x^j$ the
		Lie bracket is given by	$[V, W] = VW - WV = (VW^j - WV^j) \del /\del x^j$.
		We calculate:
		$[X_1,X_2] = -2X_3, [X_1,X_3] = 2X_2, [X_2,X_3] = -2X_1$.
		We will be interested in metrics on $M$ of the form:
		\beq
			\label{metric}
			g = dt^2 + f^2(t)(\omega^1)^2 + h^2(t)[(\omega^2)^2 + (\omega^3)^2 ] ,
		\eeq
		where $\{\omega^i\}_{i=1}^3$ is the dual coframe to $\{X_i\}$,
		$t \in I$, and $f$ and $h$ are smooth functions of $t$. 
		A metric $g$ on a smooth manifold $M$ is an Einstein metric, and the pair $(M,g)$ an Einstein manifold,
		if for all $p \in M$ and for all $X,Y \in T_p M,$
		\beqstar
			\text{Ric}(X,Y) = \lambda \cdot g(X,Y)
		\eeqstar
		for some $\lambda \in \mathbb{R}$, the Einstein constant. 
	    The Ricci curvature tensor Ric is given by the first contraction of the Riemann curvature tensor $R$, or 
		alternatively as a sum of sectional curvatures:
		\beqstar
			\text{Ric}(X, Y) = \sum_i{g(R(e_i, X)Y, e_i)},
		\eeqstar
		where $\{e_i\}_{i=0}^{n-1}$ is an orthonormal basis for $T_pM$.

		We set $e_0 = \del/\del t, e_1 = f^{-1}X_1, e_2 = h^{-1}X_2, e_3 = h^{-1}X_3$,
		and then $\{ e_i \}_{i=0}^3$ is a global, orthonormal frame for $(M,g)$.
		Note that $[\del /\del t, X_i] = 0$ for $i = 1,2,3$.
		It follows from this and the product rule that 
		\begin{align*}
			[e_0,e_1]& = -\frac{f'}{f}e_1&   [e_0,e_2] =& -\frac{h'}{h}e_2&  [e_0,e_3]& = -\frac{h'}{h}e_3 \\
			[e_2,e_3]& = -2\frac{f}{h^2}e_1&  [e_1,e_2] =& -\frac{2}{f}e_3&   [e_1,e_3]& = \frac{2}{f}e_2.		
		\end{align*}
		We are now in a position to calculate the Riemann curvature tensor $R$, given by
		\beqstar
			R(X,Y)Z = \nabla_X\nabla_Y Z - \nabla_Y\nabla_X Z - \nabla_{[X,Y]} Z, \label{curvature}
		\eeqstar
		via the Koszul formula (\cite{Kuhnel} p.214):
		\beqstar
			\bra \nabla_X Y,Z \ket = \frac{1}{2} \{ X \bra Y,Z \ket + Y\bra X,Z \ket - Z\bra X,Y \ket - 
			\bra Y,[X,Z] \ket - \bra X,[Y,Z] \ket - \bra Z,[Y,X] \ket \}, \label{Koszul}
		\eeqstar
		where here and throughout $\bra X,Y \ket = g(X,Y)$ for any tangent vectors $X$ and $Y$, and $\nabla$ is the 
		Levi-Civita connection. 
		Working in the basis $\{ e_i \}$ the matrix $(g_{ij})$ becomes the identity. The form of $g$ guarantees that all the
		off-diagonal terms of $(\text{Ric}_{ij})$ vanish, so for $g$ to be an Einstein metric it is necessary and
		sufficient that each
		diagonal element of the Ricci matrix satisfy
		$\text{Ric}_{ii} = \text{Ric}(e_i,e_i) = \lambda$. Since $X_2$ and $X_3$ are 
		treated identically in terms of the metric, $\text{Ric}_{22} = \text{Ric}_{33}$, and we have three total
		diagonal elements to consider.
		
		
		We begin with		
		\beqstar
			\text{Ric}_{00} = \bra R(e_0,e_0)e_0, e_0 \ket + \bra R(e_1,e_0)e_0, e_1 \ket + \bra R(e_2,e_0)e_0, e_2 \ket
				+ \bra R(e_3,e_0)e_0, e_3 \ket.
		\eeqstar
		Since $R(X,X)X = 0$ for any $X$, the first term is zero. For the next three terms we get
		$\bra R(e_1,e_0)e_0, e_1 \ket = -f''/f$, and
		$\bra R(e_2,e_0)e_0, e_2 \ket = \bra R(e_3,e_0)e_0, e_3 \ket = -h''/h$, which gives 
		\beqstar
			\text{Ric}_{00} = -\frac{f''}{f} - 2\frac{h''}{h}.
		\eeqstar
		Next we have 
		\beqstar
			\text{Ric}_{11} = \bra R(e_0,e_1)e_1, e_0 \ket + \bra R(e_2,e_1)e_1, e_2 \ket
				+ \bra R(e_3,e_1)e_1, e_3 \ket.
		\eeqstar
		We calculate $\bra R(e_0,e_1)e_1, e_0 \ket = -f''/f$, and
		$\bra R(e_2,e_1)e_1, e_2 \ket = \bra R(e_3,e_1)e_1, e_3 \ket = -(f'h'/fh) + (f^2/h^4)$, so
		\beqstar
			\text{Ric}_{11} = -\frac{f''}{f} - 2\frac{f'h'}{fh} + 2\frac{f^2}{h^4}.
		\eeqstar
		Finally,
		\beqstar
			\text{Ric}_{22} = \text{Ric}_{33} = \bra R(e_0,e_2)e_2, e_0 \ket + \bra R(e_1,e_2)e_2, e_1 \ket
				+ \bra R(e_3,e_2)e_2, e_3 \ket.
		\eeqstar
		We calculate $\bra R(e_0,e_2)e_2, e_0 \ket = -h''/h$, $\bra R(e_1,e_2)e_2, e_1 \ket = -(f'h'/fh) + (f^2/h^4)$, 
		and $\bra R(e_3,e_2)e_2, e_3 \ket = -(h'/h)^2 - (3f^2/h^4) + (4/h^2)$, yielding
		\beqstar
			\text{Ric}_{22} = \text{Ric}_{33} = 
				-\frac{h''}{h} - \frac{f'h'}{fh} - \frac{h'^2}{h^2} + \frac{4}{h^2} - 2\frac{f^2}{h^4} = \lambda.
		\eeqstar
		Thus the metric given by \eqref{metric} is Einstein if and only if $f$ and $h$ satisfy the following three ordinary 
		differential equations (see also \cite{Besse} p. 274):
		
		\begin{align}
			-\frac{f''}{f}& - 2\frac{h''}{h} = \lambda, \label{B1}\\
			-\frac{f''}{f}& - 2\frac{f'h'}{fh} + 2\frac{f^2}{h^4} = \lambda, \label{B2}\\
			-\frac{h''}{h}& - \frac{f'h'}{fh} - \frac{h'^2}{h^2} + \frac{4}{h^2} - 2\frac{f^2}{h^4} = \lambda. \label{B3}
		\end{align}
		So the classification of Einstein manifolds with metrics of the type given by \eqref{metric} is equivalent to 
		a classification of solutions to \eqref{B1}, \eqref{B2}, and \eqref{B3}. We shall be interested in
		solutions of this system on a maximal interval $I$ such that both $f$ and $h$ are positive.
		Because a singularity may occur at the endpoint(s) of $I$, solutions to the above system do not, in general, yield 
		geodesically complete manifolds.  However, if suitable
		boundary conditions are satisfied at the endpoint, the singularity can be removed by a change of coordinates. 
		
	\subsection*{Smoothness and completeness}
		
	At the endpoint $t = 0$ we require all of $S^3$ to collapse to a point, so $f(0) = h(0) = 0$. This is 
		analogous to the origin in polar or spherical coordinates. To avoid a conical singularity 
		we require $f'(0) = 1$, and $h'(0) = 1$ or 0.
		In cases where $f$ and $h$ are positive for all $t > 0$, there are no further restrictions on the 
		functions. If $f$ has positive roots we call the least root $l$ 
		and say the fiber collapses at $l$. In this case we require $f'(l) = -1$ to have a smooth metric. Note that if
		$h = 0$ at some point then we must have $f = 0$ or we do not have a smooth 4-manifold.
		
		A Riemannian manifold $(M,g)$ is \textit{geodesically complete} if every arc-length parameterized 
		geodesic with parameter $\tau$ is defined for all 
		$\tau \in \mathbb{R}$. The only concern in our case regards the $t$-lines of $M$, which are always geodesics 
		with arc-length parameter 
		$t$.  $(M,g)$ admits a smooth completion as long as the metric functions collapse smoothly at
		$t = 0$ and possibly at $t = l$ as described above, since we can add points to $M$ at the
		end points of $I$. In this case we'll say that $g$ or $M$ is smooth and complete.
The geometry is incomplete if, for example, $f \rightarrow \infty$ as 
		$t \rightarrow 0$. It may happen that the above conditions are satisfied except that 
		$f'(0) = n$ or $f'(0) = -f'(l) = n$ for some positive integer $n > 1$, in which case 
		we can complete $M$ to an Einstein orbifold $I \times S^3/\mathbb{Z}_n$, which is
		smooth except at the singularities at the end points, which correspond to fixed points
		of a $\mathbb{Z}_n$ action. In this case we can say $g$ or $M$ is complete but not smooth.

\subsection*{Solving the equations}	
	To solve the above system, we follow \cite{Besse} and examine the equations 
		$\eqref{B4} = \eqref{B1} - \eqref{B2}$ and  $\eqref{B5} = \eqref{B2} - \eqref{B1} + 2\cdot \eqref{B3}$:
		\begin{align}
			\frac{h''}{h}& - \frac{f'h'}{fh} + \frac{f^2}{h^4} = 0, \label{B4}\\
			-2\frac{f'h'}{fh}& - \frac{h'^2}{h^2} - \frac{f^2}{h^4} + \frac{4}{h^2} = \lambda. \label{B5}
		\end{align}
		It turns out that \eqref{B4} can be integrated to yield $f$ in terms of $h$: 
		\begin{equation}
			f = \frac{|hh'|}{\sqrt{1 - ah^2}} \label{fsolution},
		\end{equation}
		where $a \in \mathbb{R}$ is a suitable constant. Now if $a > 0$, then without loss of generality we may assume 
		$a = 1$ by the following argument. Given a metric
		that is a solution when $a > 0$, one can define the following rescaling 
		\begin{align*}
			dt& \leftrightarrow a^{-1/2}dk& f& \leftrightarrow a^{-1/2}F& h& \leftrightarrow a^{-1/2}H.
		\end{align*}
		This leaves the metric unchanged in form 
		since $G = dk^2 + F^2(\omega^1)^2 + H^2[(\omega^2)^2 + (\omega^3)^2 ] = ag$,
		and $F$ and $H$ clearly satisfy \eqref{fsolution} since $dh/dt = dH/dk$. 
		Similarly if $a < 0$ we may assume $a = -1$ and so our solutions break into three cases,
		$a = 1, -1,$ or $0$.
		
	We now introduce a new independent variable $r$ such that 
\begin{equation}
\label{cov}
		\frac{dr}{dt} = \frac{f}{h^2}.
\end{equation}
Note that this implies $dt/dr = h^2/f$ since $r$ is evidently an increasing function on $I$ and thus
		one-to-one. Equation \eqref{fsolution} then becomes autonomous in $h$ as
		\begin{align}
			\label{hsolution}
			\dot{h}& = \pm h\sqrt{1-ah^2} & \Rightarrow  r =& \pm \int{\frac{dh}{h \sqrt{1-ah^2}}},
		\end{align}
		where here and througout a dot denotes differentiation with respect to $r$.
		For the case $a = 0$ the solution is clearly $h = e^{\pm r}$, and 
		we obtain $h = \sech{r}$ for $a = 1$ and $h = \csch{r}$ for $a = -1$.
		Note that we have suppressed the constants of integration because they may be removed by a translation
		of $r$. We now proceed to find solutions for $f$ and discuss the resulting metrics.

\section{Solutions when $a = 0$}
	For this case we have $h = \dot{h} = e^r$. Substituting into \eqref{B5} yields
		\beq
			f = \sqrt{(-\lambda /6)e^{4r} + e^{2r} + Ce^{-2r}},
		\eeq
		where $C \in \mathbb{R}$ is a constant of integration. The value of the Einstein constant $\lambda$ is only
		geometrically significant up to sign. In fact, we may assume $|\lambda| = 6$ or 0, since we can rescale the metric
		as with the constant $a$, which is here zero.
		For convenience, we similarly consider cases based on the sign of $C$ and so have 9 total subcases to consider.
		
	\subsection{$\lambda = 0, C = 0$}
	
	In this case $f = h = e^r$ and 
		\beqstar
			f' = \frac{df}{dt} = \frac{df}{dr} \frac{dr}{dt} = e^r \frac{e^r}{e^{2r}} = 1 = h',
		\eeqstar
	so up to a translation in $t$ the metric is
	$g = dt^2 + t^2 [(\omega^1)^2 + (\omega^2)^2 + (\omega^3)^2]$ and $I = (0,\infty)$. $g$ is the Euclidean metric on 
	$\mathbb{R}^4$ in polar coordinates with $t$ as the distance function from the origin. The singularity at $t = 0$ 
	can be removed by changing to rectangular coordinates. 
	
	\subsection{$\lambda = 6, C = 0$}
	
	Let $s = e^r$ (we use this substitution throughout), then we have $f = \sqrt{s^2 - s^4} = s\sqrt{1 - s^2}$,
	$\frac{dt}{ds} = s/f = (1 - s^2)^{-1/2}$, which has solution 
	$s = \sin{t}$ and thus $f = \sin{t}\cos{t} = (1/2)\sin{2t}$.
	Thus
		\beq
			g = dt^2 + (1/4)\sin^2{2t}(\omega^1)^2 + \sin^2{t}[(\omega^2)^2 + (\omega^3)^2].
		\eeq
	Here $t \in (0, \pi/2)$ and
	we have an explicit smooth metric on all of $M$, which is in fact the Fubini-Study metric on $\mathbb{C} P^2$. The 
	singularity at $t = 0$ is of the same type as in section 2.1, while at $t = \pi/2$ each of the Hopf fibres is 
	collapsed to a point, but this singularity can be similarly removed by a change of coordinates. 
	
	\subsection{$\lambda = -6, C = 0$}

	In this similar case $f = s\sqrt{1 + s^2}$, and $t = \sinh^{-1}{s}$ so $f = \sinh{t}\cosh{t}$.
	This gives 
		\beq
			g = dt^2 + (1/4)\sinh^2{2t}(\omega^1)^2 + \sinh^2{t}[(\omega^2)^2 + (\omega^3)^2],
		\eeq
	and $I = (0,\infty)$. $M$ is a non-compact hyperbolic dual to the preceeding case $\mathbb{C} P^2$. The singularity
	is similarly removed at $t = 0$ and $g$ is the complex hyperbolic metric (see \cite{Kob}).
		
	\subsection{$\lambda = 0, C > 0$}
	
	We get $g = s^4(s^4 + C)^{-1}ds^2 + (s^2 + Cs^{-2})(\omega^1)^2 + s^2 [(\omega^2)^2 + (\omega^3)^2]$ and
	$t = \int_0^s{du/\sqrt{1 + Cu^{-4}}}$ with $t \in  (0,\infty)$. Note that $f \rightarrow \infty$ as 
	$s \rightarrow 0$, and since 
	$t \rightarrow 0$ as $s \rightarrow 0$, the fiber does not collapse at $t = 0$ and the metric is not complete.
	
	\subsection{$\lambda = 0, C < 0$}
	
	In this case we have the same metric as above except with negative $C$. To check the behavior at the
	endpoints, first let $D = (-C)^{1/4}$, then 
	$dt/ds = s^2(s^4 - D^4)^{-1/2}$, and
		\beqstar
			t(s) = \int_D^s{\frac{u^2 du}{\sqrt{u^4 - D^4}}} \Rightarrow f'\bigr\rvert_{t=0} = 
				\frac{df}{ds}\Bigr/\frac{dt}{ds}\biggr\vert_{s = D},
		\eeqstar
	where $t \in (0,\infty)$ and $s \in (D, \infty)$. We have
	$f' = \frac{df}{ds}\bigr/\frac{dt}{ds} = 1 + (D/s)^4$, so $f'(t=0) = 2$, 
	$h'(t=0) = 0$
	and $h(t=0) = D$. After removing the singularity at $t = 0$, we obtain a Ricci-flat Einstein metric on
	the tangent bundle of the 2-sphere $TS^2$ known as the Eguchi-Hanson metric (see \cite{Gibb}). 
	The limit metric as $D \rightarrow 0$ is
	the Euclidean metric on $\mathbb{R}^4$. Also note that for different $D > 0$, the metrics are homothetic.
	
	\subsection{$\lambda = -6, C > 0$}
	
	Here $g = s^4(s^4 + s^6 + C)^{-1}ds^2 + (s^2 + s^4 + Cs^{-2})(\omega^1)^2 + 
	s^2[(\omega^2)^2 + (\omega^3)^2]$, and $t \in (0,\infty)$, but here
	$f \rightarrow \infty$ as $t \rightarrow 0$, so $g$ is not complete.
	
	\subsection{$\lambda = -6, C < 0$}
	
	Let $D = -C$, then $g = s^4(s^4 + s^6 - D)^{-1}ds^2 + (s^2 + s^4 - Ds^{-2})(\omega^1)^2 + 
	s^2[(\omega^2)^2 + (\omega^3)^2]$. It is clear that $dt/ds = s^2(s^4 + s^6 - D)^{-1/2}$ has a single positive root, 
	call it $z$, so $I = (0,\infty)$. 
	Here $f' = \frac{df}{ds}\bigr/\frac{dt}{ds} = 1 + 2s^2 + Ds^{-4}$, and for a smooth metric we require 
	$(i.)$ $1 + 2z^2 + Dz^{-4} = n$, for an integer $n \geq 0$.
	Since we are looking for complete metrics, we may assume that $f$ vanishes at the root, so 
	$(ii.)$ $z^2 + z^4 - Dz^{-2} = 0$. Taking $(i.)$ + $z^{-2}(ii.)$ we obtain $3z^2 = n - 2$.
	In terms of $z$, $D = z^4 + z^6$, and since $h$ behaves properly as $t \rightarrow 0$,
	we have a complete Einstein metric when $D = (1/3^2)(n - 2)^2 + (1/3^3)(n-2)^3$
	over a family of manifolds $I \times S^3 / \mathbb{Z}_n$ for integers $n \geq 3$. 
	
	\subsection{$\lambda = 6, C > 0$}
	We have $g = s^4(s^4 - s^6 + C)^{-1}ds^2 + (s^2 - s^4 + Cs^{-2})(\omega^1)^2 + s^2[(\omega^2)^2 + (\omega^3)^2]$ 
	and $I = (0,\infty)$, but here as in cases 2.4 and 2.6,
	$f \rightarrow \infty$ as $t \rightarrow  0$, so $g$ cannot be complete.
	
	\subsection{$\lambda = 6, C < 0$}
	Here again let $D = -C$, and we have the same metric as above with $t(s) = \int_{z_1}^{s}{u^2du/\sqrt{u^4 - u^6 - D}}$, 
	and $I = (0, l = t(z_2))$, $s \in (z_1,z_2)$ and $z_1$ and $z_2$ are the two positive roots of $f$.
	Clearly they are both roots of $s^4 - s^6 - D$ and hence $D = z_1^4 - z_1^6 = z_2^4 - z_2^6$. For a smooth metric we
	require $f'$ to be an integer at the end points, and it follows that 
	$2 - 4z_1^2 + 2Dz_1^{-4} = n = -(2 - 4z_2^2 + 2Dz_2^{-4})$. Substituting for $D$ we obtain $z_1^2 = (4-n)/6$ and
	$z_2^2 = (4+n)/6$, so $n \leq 3$. Substituting back into the equation for $D$ we obtain 
	$6(4-n)^2 - (4-n)^3 = 6(4+n)^2 - (4+n)^3$, which is not solved for $n = 1, 2$ or $3$, and hence there are no complete
	Einstein metrics in this case.

\section{Smoothness and completeness criteria when $a \neq 0$}
	\label{Crit}
    We consider the criteria for $g$ to be complete when $a = \pm 1$. 
	From \eqref{hsolution} we know
    $h = \sech r$ for $a = 1$ and $h = \csch r$ for $a = -1$. In both cases, substituting into \eqref{B5}
    gives us the following metric
    \begin{equation}
		\label{nonzeroametric}
        g = \frac{192s^4}{(s^2+a)^2G(s^2)}ds^2 + \frac{G(s^2)}{12s^2(s^2+a)^2}(\omega^1)^2
        + \frac{4s^2}{(s^2+a)^2}[(\omega^2)^2+(\omega^3)^2],
    \end{equation}
    where $s = e^r$ and
    \begin{equation}
		\label{Gequation}
        G(x) = (24+3C-8a\lambda)(x^4+2ax^3)+48x^2+(24-3C-8a\lambda)(2ax+1).
    \end{equation}
    To determine whether $M$ is smooth and complete, we must analyze the behavior of $g$ at the endpoints of $I$.

    \subsection*{Behavior of $g$ at the endpoints}
        We observe that $f$ and $h$ are zero or infinite only if $s = 0$, $s = 1$, $s \rightarrow \infty$, or
        $G(s^2) = 0$. When $a = -1$, $s \rightarrow 1$ corresponds to $t \rightarrow \infty$. To show this when
        $G(0) > 0$, we shall show that
        \begin{equation*}
            \int_{1-\varepsilon}^1 \frac{8\sqrt{3}s^2}{(s^2+a)\sqrt{G(s^2)}}ds
        \end{equation*}
        diverges for all $\varepsilon > 0$. Consider the integrand, which is equal to $dt/ds$, as the product of
        $(s-1)^{-1}$ and a function that is continuous on $[1-\varepsilon,1]$. By the properties of continuous functions,
        $dt/ds < K (s-1)^{-1}$ for all $s \in [1-\varepsilon,1)$ and some $K > 0$. Since the integral of $(s-1)^{-1}$
        over $[1-\varepsilon,1)$ diverges, we conclude that the above integral diverges. This sort of argument also shows
        that $s = 0$ and $s \rightarrow \infty$ correspond to $t = 0$ or $t = l$, and that any root of $G(s^2)$ corresponds
        to $t = 0$ or $t = l$ if and only if it has a root of multiplicity of $1$.

        In order for $M$ to be complete we require $f = 0$ at $t = 0$ and $t = l$, so if $G(0) > 0$ $g$ is not complete
        since $f \rightarrow \infty$ as $s \rightarrow 0$. Similarly, $g$ is not complete when $24-3C-8a\lambda > 0$
        since $f \rightarrow \infty$ as $s \rightarrow \infty$. Otherwise $f(t) = 0$ and $h(t)$ is $0$ or finite at
        $t = 0$ and $t = l$, so $g$ is complete provided $M$ is a smooth manifold. We have $df/dt = dh/dt = 1$ at
        $s = 0$ when $G(0) = 0$, so $M$ is smooth at $s = 0$. Furthermore, when $24-3C-8a\lambda = 0$, $M$ is smooth for
        $s$ near $\infty$ since $df/dt$ and $dh/dt$ tend to $-1$ as $s \rightarrow \infty$. Thus we must determine if
        $M$ is smooth at roots of $G(s^2)$.

    \subsection*{Smoothness at one root of $G$}
        Suppose $I = (0,\infty)$ and let $z > 0$ satisfy $G(z^2) = 0$. Since $dh/dt =0$ at $s = z$, to show that $M$
        is smooth we must show that for some integer $n$.
        \begin{equation}
            \label{fderivative}
            \frac{df}{dt} = \frac{(24+3C-8a\lambda)(2z^6+3az^4)+48z^2+a(24-3C-8a\lambda)}{24z^2} = n.
        \end{equation}
        Since $G(z^2) = 0$, whenever $z \neq 1$ we have
        \begin{equation}
            \label{Cvalue1}
            C = \frac{(24-8a\lambda)(z^8+2az^6+2az^2+1)+48z^4}{-3(z^8+2az^6-2az^2-1)}.
        \end{equation}
        By substituting for $C$ in \eqref{fderivative} and solving for $\lambda$, we obtain
        \begin{equation}
            \label{lambda1}
            \lambda = \frac{(2+n)z^4+4az^2+(2-n)}{2z^2}.
        \end{equation}
        Thus for any given choice of $z$ and $n$, there exist values for $C$ and $\lambda$ such that the manifold is
        smooth at $s = z$.

    \subsection*{Smoothness at two roots of $G$}
        Suppose $I = (0,l)$ and let $z_1$ and $z_2$ be roots of $G(s^2)$. We would like to show that $df/dt = n$ at 
		$s = z_1$ and $df/dt = -n$ at $s = z_2$ for some positive integer $n$.
        Thus we shall consider the following system of equations, where $i = 1$, $2$,
        \begin{eqnarray*}
            G(z_i^2) &=& 0,\\
            \frac{df}{dt} &=& (-1)^{i+1} n ~\mbox{at}~ z_i.
        \end{eqnarray*}
        As before, we can now solve for $C$ and $\lambda$ in terms of $z_1$ and $n$. Whenever
        $z_1 \neq 1$, we have
        \begin{eqnarray}
            \label{Cvalue2}
            C &=& \frac{(24-8a\lambda)(z_1^8+2az_1^6+2az_1^2+1)+48z_1^4}{-3(z_1^8+2az_1^6-2az_1^2-1)},\\
            \label{lambda2}
            \lambda &=& \frac{(2+n)z_1^4+4az_1^2+(2-n)}{2z_1^2}.
        \end{eqnarray}

        Now consider $df/dt = -n$ at $z_2$. Observe that if we fix $\lambda$ and solve for $z_1$, we have have a
        quadratic equation in $z_1^2$. Thus the only non-negative solutions to this equation are $z_1$ and
        $\sqrt{\frac{2-n}{2+n}} z_1$. By replacing $n$ with $-n$ in \eqref{lambda2}, we observe that
        $z_2^{-1}$ is a solution to \eqref{lambda2} for a fixed $\lambda$. If $z_1 = z_2^{-1}$, then by
        substituting in for $z_1$ into \eqref{Cvalue2}, we obtain $C = -C$, thus $C = 0$. Otherwise,
        $z_2 = \sqrt{\frac{2+n}{2-n}} z_1^{-1}$, which we will assume for the remainder of this discussion.

        We must choose $z_1$ and $n$ such that $C$ is the same if we replace $z_1$ with $z_2$ in
        \eqref{Cvalue2}. Let $C_1$ denote the value of $C$ given by \eqref{Cvalue2} and let $C_2$ denote the
        value of $C$ when $z_1$ is replaced by $z_2$. We have
        \begin{equation*}
            C_1 - C_2 = \frac{32n^3 z_1^4}{3(2a(1+az_1^2)^2+n(1-az_1^4)^2)}.
        \end{equation*}
        Thus $C_1 = C_2$ if and only if $n = 0$ or $z_1 = 0$. Since $n = 0$ provides us with a root of multiplicity
        greater than one and $z > 0$ by assumption, the manifold is smooth only when $z_1 = 1$ or $C = 0$.

\section{Solutions when $a = 1$}
    For this case we have the following metric $g$ on $M$,
    \begin{equation}
        g = \frac{192s^4}{(s^2+1)^2G(s^2)}ds^2 + \frac{G(s^2)}{12s^2(s^2+1)^2}(\omega^1)^2
        + \frac{4s^2}{(s^2+1)^2}[(\omega^2)^2+(\omega^3)^2],
    \end{equation}
    where $G(x)$ is given by \eqref{Gequation} with $a = 1$.
    We may assume $C \geq 0$ since if $C < 0$, we can make the change of variable $u \leftrightarrow s^{-1}$ to
    produce a metric where $C$ is replaced by $-C$ and that is otherwise identical to $g$. We have the following three
    subcases, determined by the behavior of $g$ at the endpoints of $I$.

    \subsection{$24+3C-8\lambda > 0$}
        If $24-3C-8\lambda \geq 0$, $G(s^2) > 0$ for all $s > 0$. If $24-3C-8\lambda < 0$, by Descartes' Rule of Signs,
        $G(s^2)$ has exactly one positive root, $z$ and $G(s^2) > 0$ on $(z,\infty)$. In both cases
        $f \rightarrow \infty$ as $s \rightarrow \infty$ and thus $g$ is not complete. Observe that $a =1$ and
        $\lambda \leq 0$ only in this case, so all complete solutions for $a = 1$ occur when $\lambda > 0$.

    \subsection{$24+3C-8\lambda = 0$}
        By the quadratic formula, $G(s^2)$ has exactly one positive root, $z$, which satisfies
        $z^2 = \frac{1}{8}C+\frac{1}{8}\sqrt{C^2+8C}$. Thus $G(s^2) > 0$ on $(z,\infty)$ and $I = (0,l)$ where
        $f(0) = f(l) = 0$, so $M$ is complete and we proceed to check smoothness. 
		Since $df/dt = 1$ at $s = 0$, we require
        $df/dt = -1$ at $s = z$. Solving produces exactly one solution at $C = 0$ and $\lambda = 3$. In this case we have
        \begin{equation*}
            \frac{dt}{ds} = \frac{2}{s^2+1} \Rightarrow s = \tan(t/2) \Rightarrow f = h = \sin(t).
        \end{equation*}
        Thus we get the explicit solution $g = dt^2+\sin^2(t)[(\omega^1)^2+(\omega^2)^2 +(\omega^3)^2]$, where
        $t \in (0,\pi)$. $(M, g)$ is the unit 4-sphere.

    \subsection{$24+3C-8\lambda < 0$}
    \label{Page's example}
        By Descartes' Rule of Signs, $G(s^2)$ has at most two positive roots, $z_1$ and $z_2$. To determine when $G$
        has two roots, we shall first determine the boundary of the regions in the $C$-$\lambda$-plane where $G$ has
        no roots and where $G$ has two roots. This boundary occurs when $G(x_0) = G'(x_0) = 0$ for some $x_0 > 0$.
        Solving gives us the following relationship,
        \begin{equation*}
            C = \frac{8(3-\lambda) \left[8-(\lambda-2\sqrt{\lambda^2-4\lambda})^3 \right]}
            {24+3(\lambda-2\sqrt{\lambda^2-4\lambda})^3}.
        \end{equation*}
        Let $C_0$ denote this value for $C$. Since $x_0 < 1$, $G(x_0)$ decreases as $C$ increases. Thus when $C < C_0$,
        $G(x_0) > 0$ and $G$ has exactly two positive roots. When $C \leq C_0$, $G(x_0) \leq 0$, so $G(x) \leq 0$ for all
        $x > 0$ since $G$ has a maximum at $x = x_0$. In cases where two roots occur, $g$ is defined for
        $s \in (z_1,z_2)$. Thus $I = (0,l)$ and from section \ref{Crit} we know 
		the manifold is smooth only if $z_1 = 1$ or $C = 0$.
        The case where $z_1 = 1$ yields no smooth manifolds. However, when $C = 0$, we obtain a solution with
        $df/dt = 1$ at $s = z_1$ when $z_1^{-1} = z_2 = \frac{1}{2}\sqrt{y}+\frac{1}{2}\sqrt{-y+8/\sqrt{y}}$, where
        $y = 2(\sqrt[3]{1+\sqrt{2}}+\sqrt[3]{1-\sqrt{2}})$. The root $z_2$ was calculated using the formula for solving
		general cubic equations.
		In this case, $\lambda$ is given by \eqref{lambda2}
        with $n = 1$. 
        The manifold $(M,g)$ is a non-trivial $S^2$ bundle over $S^2$, topologically
        $\mathbb{C}P^2 \# \overline{\mathbb{C}P^2}$. This is the example originally given by D. Page \cite{Page}.

\section{Solutions when $a = -1$}
    For this case we have
    \begin{equation}
        g = \frac{192s^4}{(s^2-1)^2G(s^2)}ds^2 + \frac{G(s^2)}{12s^2(s^2-1)^2}(\omega^1)^2
        + \frac{4s^2}{(s^2-1)^2}[(\omega^2)^2+(\omega^3)^2],
    \end{equation}
    and $G(x)$ is given by \eqref{Gequation} with $a = -1$.
    We may assume $C \geq 0$ for the same reason as when $a = 1$. We have the following
    subcases.

    \subsection*{Some families of solutions when $\lambda \leq 0$}
    We begin by considering three families of complete solutions.

    \subsection{$24+3C+8\lambda < 0$}
        By Descartes' Rule of Signs , one readily determines that $G$ has at most two positive roots. Since
        $G(0) < 0$, $G(1) > 0$, and $G(s^2) \rightarrow -\infty$ as $x \rightarrow \infty$, by the intermediate 
        value theorem
        $G(s^2)$ has exactly two positive roots, $z_1$ and $z_2$, such that $z_1 < 1 < z_2$. Thus $f > 0$ for
        $s \in (z_1,1) \cup (1,z_2)$, and when one restricts $s$ to either $(z_1,1)$ or $(1,z_2)$, $I = (0,\infty)$.
        Thus $g$ is complete provided it is smooth at $t = 0$.

    \subsection{$24+3C+8\lambda = 0$}
        Since $G$ is a quadratic polynomial, we readily determine that $G(s^2)$ has exactly one positive root,
        namely
        \begin{equation*}
            z^2 = \frac{-C+\sqrt{C^2+8C}}{8} < 1.
        \end{equation*}
        The metric $g$ is defined for any $s \in (z,1) \cup (1,\infty)$. When $s$ is restricted to $(z,1)$,
        $I = (0,\infty)$ and thus $g$ is complete provided it is smooth at $t = 0$. When $s$ is
        restricted to $(1,\infty)$, $I = (0,\infty)$ and $g$ is complete.

    \subsection{$3C > |24+8\lambda|$}
        Since $G(0) < 0$ and $G(1) > 0$, $G(s^2)$ has at least one root on $(0,1)$. We know $G$
        is positive on $(1,\infty)$ since $G(1) > 0$ and $G'(x) > 0$ on $(0,1)$. Furthermore, $G$ has no more than one
        root on $(0,1)$ because if a minimum occurs at $x > 0$ such that $G(x) < 0$, then $x > 1$ since
        \begin{equation}
            (24+3C+8\lambda)(x-1)(x+1)^3 > \frac{1}{2} G'(x) - G(x) > 0.
        \end{equation}
        Thus the metric $g$ is defined for any $s \in (z,1) \cup (1,\infty)$. When $s$ is restricted to $(z,1)$,
        $I = (0,\infty)$ and thus $g$ is complete provided it is smooth at $t = 0$. When $s$ is restricted to
        $(1,\infty)$, $g$ is not complete. \\

        In the previous three sections, we concluded that $g$ is complete provided it is smooth at $t = 0$. In fact, 
		for each positive integer $n$ we have a continuous family of manifolds as $z$ ranges across $(0,1)$ and with $C$
        and $\lambda$ given by \eqref{Cvalue1} and \eqref{lambda1}. Also, given an integer $n < 0$ we obtain other
        continuous families of manifolds dependent on $z \in (1,\infty)$. We require that $z^2 \geq 1/3$ when $n = 1$
        and $z^2 \leq 3$ when $n = -1$. When $z^2 = 1/3$ and $n = 1$, $\lambda = 0$ and $C = 10$. We obtain no other
        smooth complete metrics for $\lambda = 0$ and $C > 8$.

    \subsection*{Other solutions when $\lambda \leq 0$}
    \subsection{$24-3C+8\lambda = 0$}
        Here $G$ has no positive roots, so $G(x) > 0$ on $[0,\infty)$. Thus $g$ is defined for
        $s \in (0,1) \cup (1,\infty)$. When $s$ is restricted to $(0,1)$, $I = (0,\infty)$ and $g$ is complete. Thus
        we obtain a continuous family of smooth and complete manifolds dependent on $\lambda$ and with $C = (24+8\lambda)/3$.
        Note that this family contains the manifold with $\lambda = 0$ and $C = 8$.

        We can explicitly write the metric for one manifold in this family in the case that $C = 0$ and $\lambda = -3$.
        Restricting $s$ to $(0,1)$, we obtain
        \begin{equation*}
            \frac{dt}{ds} = \frac{2}{s^2-1} \Rightarrow s = \tanh(t/2) \Rightarrow f = h = \sinh(t).
        \end{equation*}
        We get $g = dt^2+\sinh^2(t)[(\omega^1)^2+(\omega^2)^2 +(\omega^3)^2]$,
        where $t \in (0,\infty)$. We derive the same metric in the case that $s$ is restricted to $(1,\infty)$.
        $(M,g)$ is the complete hyperbolic 4-manifold with constant curvature $-1$.

        When $s$ is restricted to $(1,\infty)$ and $C > 0$, $g$ is not complete since $f \rightarrow \infty$ as
        $s \rightarrow \infty$.

    \subsection{$24-3C+8\lambda > 0$}
		When $s\neq 1$, we have
        \begin{equation*}
            G(s^2) = (24+3C+8\lambda)(s^4-s^2)^2-16\lambda s^4+(24-3C+8\lambda)(s^2-1)^2 > 0.
        \end{equation*}
        Thus the metric $g$ is defined for any $s \in (0,1) \cup (1,\infty)$. When $s$ is restricted to either $(0,1)$
        or $(1,\infty)$, $g$ is not complete.

    \subsection*{Solutions when $\lambda > 0$}
    \subsection{$3C > 24+8\lambda$}
        In this case, $G$ may have either one or three positive roots depending on the values of $C$ and $\lambda$.
        One root, $z_1$, is guaranteed and satisfies $z_1 > 1$. The other two roots, $z_2$ and $z_3$, satisfy
        $z_2 \leq z_3 < 1$. In either case $g$ is defined for $s \in (z_1, \infty)$ and $g$ is not complete when $s$
        is restricted to $(z_1, \infty)$.

        To determine when $G$ has three roots, we shall first determine the boundary of the regions in the
        $C$-$\lambda$-plane where $G$ has one root and where $G$ has three roots. This boundary occurs when
        $G(x_0) = G'(x_0) = 0$ for some $x_0 > 0$. Solving gives us the following relationship,
        \begin{equation*}
            C = \frac{8(\lambda+3) \left[8+(\lambda-2\sqrt{\lambda^2+4\lambda})^3 \right]}
            {24-3(\lambda-2\sqrt{\lambda^2+4\lambda})^3}.
        \end{equation*}
        Let $C_0$ denote this value for $C$. For any fixed $x < 1$, $G(x)$ decreases as $C$ increases. Thus when
        $C > C_0$, $G(x) < 0$ and $G$ has only one positive root. When $C \leq C_0$, $G(x_0) > 0$, so $G$ has three
        positive roots since $G$ is negative at $s = 0$ and $s = 1$.

        When three distinct roots occur, $g$ is also defined for $s \in (z_2,z_3)$. Here $I = (0,l)$ and $(M,g)$ is
        not smooth.

    \subsection{$24-3C+8\lambda = 0$}
        We have $G(0) = 0$ and $G(1) < 0$. We find two positive roots, $z_1$ and $z_2$,
        such that $z_1 < 1 < z_2$ and
        \begin{equation}
            \label{zformula}
            z_i = 1 + (-1)^i \sqrt{1-\frac{8}{C}}, ~\mbox{where}~ i = 1,2.
        \end{equation}
        Thus $g$ is defined for $s \in (0,z_1) \cup (z_2,\infty)$. If we restrict $s$ to $(0,z_1)$, then $I = (0,l)$. 
		Since at $s = 0$, $df/dt = dh/dt = -1$, we must find a $z_1$ satisfying \eqref{zformula} and
        $df/dt = \frac{1}{4}Cz^2 = n$ at $s = z_1$. This has no solutions. If
        we restrict $s$ to $(z_2, \infty)$, we also find no solutions.

    \subsection{$24-3C+8\lambda > 0$}
        Since $G(0) > 0$, $G(1) < 0$, and $G(s^2) \rightarrow \infty$ as $s \rightarrow \infty$,
        $G$ has at most two positive roots. Also, $G''$ has positive roots
        only on $(0,1)$, so $G$ has exactly one root greater than $1$. Furthermore, $G$ has no more than one root on
        $(0,1)$ by the argument used in the case where $3C > |24+8\lambda|$ and $\lambda \leq 0$. Thus $g$ is defined for
        $s \in (0,z_1) \cup (z_2,\infty)$, where $z_1$ and $z_2$ are roots of $G(s^2)$. When $s$ is restricted to either
        $(0,z_1)$ or $(z_2,\infty)$, $g$ is not a complete metric. 
        
\section{Conclusion}

	We have been able to completely solve the problem
	of finding Einstein metrics on $M$, the cylinder over the 3-sphere, of the form given by \eqref{metric}. We
	have demonstrated that $M$ supports a wide variety of Einstein metrics and provided explicit formulas wherever
	possible. While writing this paper we found the recent work \cite{Lu et al.} which also
	contains new Einstein metrics, though we are not sure if they are identical with ours.
 A thorough investigation of the geometry of these various solutions would be interesting and 
	remains for future work.


\end{document}